 \documentclass[11pt]{amsart}
\usepackage{epsfig,graphicx}

\newcommand{\mcg}{\Gamma _{g,r}^n}

\newtheorem{theorem}{Theorem}[section]
\newtheorem{thm}[theorem]{Theorem}
\newtheorem{lem}[theorem]{Lemma}
\newtheorem{cor}[theorem]{Corollary}

\newtheorem{rem}{Remark}[section]

\begin{document}

 \title[On cofinite subgroups of mapping class groups]
  {On cofinite subgroups of mapping class groups }

 \author{Mustafa Korkmaz}



 \address{Department of Mathematics, Middle East Technical University,
 06531 Ankara, Turkey}
 \email{korkmaz@arf.math.metu.edu.tr}

 \def\R{{\mathbb{R}}}
 \def\Z{{\mathbb{Z}}}
 \def\mcg{{\rm Mod}}
 \def\pmcg{{\rm PMod}}

\maketitle

 \begin{abstract}
 For any positive integer $n$, we exhibit a cofinite subgroup $\Gamma_n$ of the
 mapping class group of a surface of genus at most two such that $\Gamma_n$ 
admits an epimorphism onto a free group of rank $n$. We conclude that 
$H^1(\Gamma_n;\Z)$ has rank at least $n$ and the dimension of the second bounded 
cohomology of each of these mapping class groups is the cardinality of the continuum.
In the case of genus two, the groups $\Gamma_n$ can be chosen not to
contain the Torelli group. Similarly for hyperelliptic mapping class groups. 
We also exhibit an automorphism of a subgroup of finite index in the mapping
class group of a sphere with four punctures (or a torus) 
such that it is not the restriction of an endomorphism of the whole group.
 \end{abstract}

\section{Introduction}
It is well known that the first homology group of the mapping class group of
a closed orientable surface of genus $g$ is trivial for $g\geq 3$ and isomorphic
to $\Z_{12}$ and $\Z_{10}$ if $g=1$ and $g=2$ respectively. It follows that
the first cohomology of this group is trivial.  N.V.~Ivanov
(Problem $2.11(A)$ in~\cite{ki}) asked whether $H^1(\Gamma;\Z)$ is trivial
for any subgroup $\Gamma$ of finite index in the mapping class group.
In the case $g\geq 3$, this question was answered affirmatively by J.~D.~McCarthy~\cite{mc} 
for subgroups $\Gamma$ containing the Torelli group, 
the subgroup of the mapping class group consisting of those mapping 
classes that act trivially on the first homology of the surface.
For arbitrary subgroup of finite index the problem is still open.
It was also shown by McCarthy~\cite{mc} and Taherkhani~\cite{t} 
that the mapping class group of  a closed 
orientable surface of genus $2$ contains subgroups of finite index with nontrivial 
first cohomology. 
All of the examples of McCarthy contain the Torelli group.
More precisely, he shows that if $r$ is an integer divisible by $2$ or $3$, then the
kernel of the action of the mapping class group on the mod~$r$ homology of the surface 
has nontrivial first cohomology.  It is not clear whether the examples of Taherkhani
contain the Torelli group, because his calculations are carried out by computer.

The purpose of this paper is to give an elementary construction of a sequence 
$\Gamma_n$ of subgroups of finite index 
in the mapping class group of an orientable surface of genus at most $2$ and
 the hyperelliptic mapping class group such that $\Gamma_n$ admits a
homomorphism onto a finitely generated free group of rank $n$. In the case of  
a closed orientable surface of genus $2$, we can choose these subgroups 
in such a way that they do not contain
the Torelli group. This shows that
for any positive integer $n$, there is a subgroup of finite index whose first cohomology 
has rank at least $n$. Another application is that the dimension of the second bounded 
cohomology of each of these mapping class groups is the cardinal of the continuum. 
The fact that they are infinite dimensional was also proved by 
Bestvina and Fujiwara~\cite{bf} by completely different arguments.

The last section is independent of the other results in the paper. 
In this section we prove that 
there is a subgroup $\Gamma$ of finite index in the 
mapping class group of a sphere with four punctures
and in that of a torus (or a torus with one puncture), 
and an automorpism $\varphi: \Gamma\to\Gamma$
such that $\varphi$ is not the restriction of any endomorphism of the whole group.
It is known that if the surface is not a sphere with four punctures
or a torus with at least two punctures, then any isomorphism 
between two subgroups of finite index in the mapping class group
is the restriction of an automorphism of the whole group (cf.~\cite{i1,ko}).  
In case a of torus with two punctures, the answer to the related obvious question
is unknown.

 \section{Definitions and preliminaries}

Let $S$ be an orientable surface of genus $g$ with $p$ marked points (=punctures) 
and with $q$ boundary components. The mapping class group $\mcg_{g,p}^q$ 
is defined to be the group of isotopy classes
of orientation preserving diffeomorphisms $S\to S$ which restrict to the identity on the boundary 
and preserve the set of punctures. The isotopies are assumed to fix the punctures and 
the points on the boundary.  If $p$ and/or $q$ is zero, then we omit it from the notion, so that 
$\mcg_{g}^q,\mcg_{g,p}$ and $\mcg_{g}$ denote $\mcg_{g,0}^q,\mcg_{g,p}^0$ and $\mcg_{g,0}^0$
respectively.

The pure mapping class group $\pmcg_{g,p}^q$ is the kernel of
of the action of $\mcg_{g,p}^q$ on the set of punctures. Since this action of $\mcg_{g,p}^q$
is transitive, the quotient of $\mcg_{g,p}^q$ by $\pmcg_{g,p}^q$ is isomorphic to 
the symmetric group on $p$ letters.

The Torelli group is the kernel of the natural map $\mcg_{g}\to Sp(2g,\Z)$ obtained from
the action of the mapping class group on the first homology of the surface $S$.
 
\begin{figure}[hbt]
    \begin{center}
          \epsfig{file=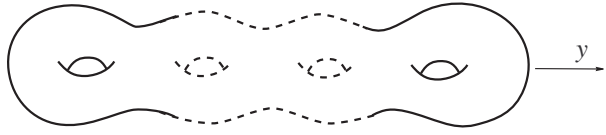,width=8.5cm}
          \caption{A closed orientable surface embedded in $\R^3$ and it is invariant under $J$.}
          \label{yuzey}
    \end{center}
 \end{figure}

Suppose now that $S$ is closed and it is embedded in the $xyz$-space as in 
Figure~\ref{yuzey} in such a way that it is invariant under the
 rotation $J(x,y,z)=(-x,y,-z)$ about the $y$-axis. Let $\jmath$ denote the isotopy class of
 $J$. The centralizer 
$$\Delta_g=\{ f\in \mcg_g \ | \ f\jmath =\jmath f \} $$
 of $\jmath$ in $\mcg_g$ is called the hyperelliptic mapping class
 group. Note that if $g=1$ or $2$, then the hyperelliptic mapping class
 group is equal to the mapping class group.

The involution $J$ has $2g+2$ fixed points. Thus the quotient of $S$
by $J$ gives a branced covering $\pi : S\to R$ branching over 
$2g+2$ points, where $R$ is the $2$-sphere. This branced
covering induces a short exact sequence
 \begin{eqnarray} \label{seq}
 1\longrightarrow \Z_2 \longrightarrow
 \Delta_g\stackrel{\pi_*}\longrightarrow \mcg_{0,2g+2} \longrightarrow 1, 
\end{eqnarray}
 where $\Z_2$ is the subgroup of $\Delta_g$ generated by the hyperelliptic involution
 $\jmath$ (cf. \cite{bh}).

 \section{Finite index subgroups with large $H^1$}

 In this section we give the construction of subgroups admitting 
homomorphisms onto free groups.

 \begin{thm}
 Suppose that $g\leq 2$. If $g=0$, suppose, in addition, that
 $p+q\geq 4$. For any positive integer $n$, there is a subgroup $\Gamma_n$ of finite
 index in $\mcg_{g,p}^q$ such that there is a epimorphism from $\Gamma_n$ onto
a free group $F_n$ of rank $n$ and $\Gamma_{n+1}\subset \Gamma_n$.  In the case of $g=2$ and $p=q=0$, 
the group $\Gamma_n$ can be chosen not to contain the Torelli group.
 \label{th:main}
 \end{thm}
 \begin{proof}
 For $n\geq 4$ it is well known that forgetting one of the punctures on 
a sphere with $n$ punctures gives rise to a short exact sequence
$$1\to F_{n-2}\to \pmcg_{0,n}\to \pmcg_{0,n-1}\to 1,$$
where $F_{n-2}$ is the fundamental group of a sphere with $n-1$
 punctures, which is a free group of rank $n-2$.
 It can easily be shown that $\pmcg_{0,3}$ is
 trivial. It follows that $\pmcg_{0,4}$ is a free group of rank $2$.

 We first prove the theorem for $n=2$. That is, we prove that there
 is a finite index subgroup $\Gamma_2$ in $\mcg_{g,p}^q$ and an epimorphism
$\Gamma_2\to F_2$.

 Suppose first that $g=0$.
 There is an epimorphism $ \pmcg_{0,p}^q\to\pmcg_{0,4}$
 obtained by gluing a disc with one puncture to each boundary
 component and then forgetting some $p+q-4$ punctures. The subgroup
 $\pmcg_{0,p}^q$ is of index $p\, !$ in $\mcg_{0,p}^q$. In this
 case, we can take $\Gamma_2$ to be the subgroup $\pmcg_{0,p}^q$.

 Suppose next that $g=1$.
 Gluing a disc along each boundary component and forgetting all
 punctures yield an epimorphism $\varphi:\mcg_{1,p}^q\to \mcg_1$.
 The group $\mcg_1$ is isomorphic to $SL(2,\Z )$. 
 The commutator subgroup $[\mcg_1,\mcg_1]$ of $\mcg_1$
 is a free group of rank $2$ and its index in $\mcg_1$ is $12$. Thus
 we can take $\Gamma_2$ to be $\varphi^{-1}([\mcg_1,\mcg_1])$.

 Suppose finally that $g=2$. Again, gluing a disc along each
 boundary component and forgetting the punctures give an epimorphism
 $\varphi :\mcg_{2,p}^q\to \mcg_2$. Note that $ \mcg_2=\Delta_2$. Consider the natural map
 $\pi_* :\mcg_{2}\to \mcg_{0,6}$ in~(\ref{seq}).  Since there is an epimorphism
 from $\pmcg_{0,6}$ onto the free group $\pmcg_{0,4}$ of rank $2$,
 we may take $\Gamma_2=\varphi^{-1} ( \pi_* ^{-1} ( \pmcg_{0,6} ) )$. The index
 of $\Gamma_2$ in $\mcg_{2,p}^q$ is $720$. 

 For $n\geq 3$, consider an epimorphism $f:\Gamma_2\to F_2$.
 Let $F_n$ be a subgroup of $F_2$ of index $n-1$. Then $F_n$ is a
 free group of rank $n$. The subgroup $\Gamma_n=f^{-1}(F_n)$
 is of finite index in $\mcg_{g,p}^q$ and the restriction of $f$
 maps $\Gamma_n$ onto $F_n$. This completes the proof of
the first assertion.

In the case $g=2$ and $p=q=0$, we can choose $\Gamma_n$ so that it does not
contain the Torelli group as follows.  Let $S$ be a closed connected oriented
 surface of genus $2$ embedded in the $xyz$-space as in Figure~\ref{yuzey}.
Let $a$ be the separating simple closed curve which is the intersection of
the $xz$-plane with $S$.  We note that $a$ passes through no fixed 
points of $J$ and $J(a)=a$. The quotient of $S$ by the action of $J$ gives rise to 
a branched covering $\pi:S\to R$, where $R$ is 
a $2$-sphere. Let us denote the image of the fixed points 
$J$ by  $P_1,P_2,\ldots,P_6$, so that we see them as punctures on $R$. 
The simple closed curve $\pi (a)$ separates the punctures on $R$ into two sets
each containing three elements. We can assume that $P_1,P_2,P_3$ is separated
from the other three punctures by $\pi (a)$. 

Now assume that $P_3$ and $P_6$ are not marked points on $R$.
Let $c$ denote the image of $\pi (a)$ on this sphere with four punctures
$P_1,P_2,P_4,P_5$. Choose an embedded arc $\delta$ on $R$ connecting 
the punctures $P_2$ and $P_4$ so that it intersects $c$ only once
and its interior does not contain any puncture. If $d$ denotes the
boundary of a regular neighborhood of $\delta$, it can easily be shown
that the Dehn twists $t_c$ and $t_d$ generate $\pmcg_{0,4}=F_2$ freely.
For $n\geq 3$, the subgroup of $F_2$ generated by 
$t_d,t_ct_dt_c^{-1}, t_c^2t_dt_c^{-2}, \ldots, t_c^{n-2}t_dt_c^{-n+2}$ and $t_c^{n-1}$
is a free group $F_n$ of rank $n$. Note that for $n\geq 4$, the element $t_c^2$
is not contained in $F_n$.

Since the curve $a$ does not contain any fixed point of $J$, the restriction 
of $\pi$ to $a$ gives an honest two-sheeted covering $a\to \pi (a)$. It is easy to see that
$\pi_* (t_a)=t_{\pi (a)}^2$, which is contained in $\pmcg_{0,6}$. 
If $\phi :\pmcg_{0,6}\to \pmcg_{0,4}$ is the epimorphism obtained by forgetting the punctures 
$P_3$ and $P_6$, then clearly we have $\phi (t_{\pi (a)})=t_c$, and 
so $\phi(\pi_* (t_a))=\phi(t_{\pi (a)}^2)=t_c^2$.

If we define $\Gamma_2$ and $\Gamma_3$ to be the subgroup $ \pi_*^{-1}(\phi ^{-1} (F_4))$
and $\Gamma_n$ to be $  \pi_*^{-1}(\phi ^{-1} (F_n))$ for $n\geq 4$, obviously there exists an
epimorphism from $\Gamma_k$ onto a free group of rank $k$ for all $k\geq 2$.  The element 
$t_a\in\mcg_2$ is contained in the Torelli group and  but not in $\Gamma_k$.

 This completes the proof of the theorem.
 \end{proof}

 \begin{rem}
 Suppose that $p+q\leq 3$. In this case the mapping class group 
$\mcg_{0,p}^q$ is 
\begin{itemize}
\item trivial if $p\leq 1$ and $q\leq 1$,
\item the cyclic group of order $2$ if $(p,q)=(2,0)$,
\item the symmetric group on $3$ letters if $(p,q)=(3,0)$,
\item $\Z$ if $(p,q)=(0,2)$ or $(2,1)$,
\item $\Z\oplus \Z$ if $(p,q)=(1,2)$, and
\item $\Z\oplus \Z\oplus\Z$ if $(p,q)=(0,3)$.
\end{itemize}
Thus non of these groups have a subgroup admitting a homomorphism
onto a free group of rank greater than $1$. 
 \end{rem}

 \begin{rem}
For each $n$, the subgroup $\Gamma_n$ of $\mcg_2$ in the above theorem
can also be chosen to contain the Torelli group.
 \end{rem}

 \bigskip

 \begin{cor}
Suppose that $g\geq 2$. For any positive integer $n$, there is a
subgroup $\Gamma_n$ of finite index in the hyperelliptic mapping class
group $\Delta_g$ such that there is an epimorphism from
$\Gamma_n$ onto a free group $F_n$ of rank $n$.
 \label{cor:main}
 \end{cor}
 \begin{proof}
The corollary follows easily from the fact that $\Delta_g$ admits an 
epimorphism onto $\mcg_{0,2g+2}$ and Theorem~\ref{th:main}.
 \end{proof}

\bigskip

The next two corollaries follow from Theorem~\ref{th:main} and Corollary~\ref{cor:main}.

 \begin{cor}
 Suppose that $g\leq 2$. If $g=0$, suppose, in addition, that $p+q\geq 4$. 
For any positive integer $n$, there is a subgroup $\Gamma_n$ of finite
index in $\mcg_{g,p}^q$ such that the rank of $H^1(\Gamma_n;\Z)$ is at
least $n$.
 \label{cor1:H^1}
 \end{cor}

 \begin{cor}
For any positive integer $n$, there is a subgroup $\Gamma_n$ of finite
index in the hyperelliptic mapping class group
$\Delta_g$ such that the rank of $H^1(\Gamma_n;\Z)$ is at
least $n$. Moreover, in the case of $g=2$ the subgroup $\Gamma_n$ of
$\mcg_2$ can be chosen so that it does not contain the Torelli group.
 \label{cor2:H^1}
 \end{cor}

 \section{The second bounded cohomology}
 In this section, we show how to deduce from Theorem~\ref{th:main}
that the dimension of the second bounded cohomology group of 
the mapping class group $\mcg_{g,p}^q$ for $g\leq 2$ 
and that of the hyperelliptic mapping class group
$\Delta_g$ is the cardinality of the continuum.

Let $G$ be a discrete group and let
 $$ C_b^k(G;\R)=\{ f:G^k\to \R \,\,|\,\, f(G^k) \mbox{ is bounded } \}.
 $$
There is a coboundary operator $\delta_b^k:C_b^k(G;\R)\to
C_b^{k+1}(G;\R)$ defined by
 \begin{eqnarray*}
 \delta_b^k (f)(x_0,\ldots, x_{k})&=& f(x_1,\ldots, x_{k})
  +\sum_{i=1}^{k} (-1)^i
 f(x_0,\ldots,x_{i-1}x_{i},\ldots, x_{k})\\
 &&+(-1)^{k+1}f(x_0,\ldots, x_{k-1}).
 \end{eqnarray*}
The cohomology of the complex $\{ C_b^k(G;\R),\delta_b^k \}$ is
called the bounded cohomology of $G$ and is denoted by
$H_b^*(G;\R)$. The space $ C_b^k(G;\R)$ is a Banach space with the
norm
 $$ ||f||=\sup \{ \,\,|f(x_1,x_2,\ldots,x_k)|\,\,:\,\, x_i\in G\}, $$
which induces a semi-norm on $H_b^k(G;\R)$. The bounded cohomology
$H_b^k(G;\R)$ is always a Banach space for $k=2$ (cf. \cite{iv2}) but it
need not be a Banach space for $k\geq 3$ (cf. \cite{s}).

The first result in the theory of bounded cohomology is that the
first bounded cohomology of any group is trivial. This is because
a bounded $1$-cochain is a bounded homomorphism $G\to \R$ and 
any such homomorphism is trivial. So the
first interesting bounded cohomology is in dimension two.

In the above definition, if we replace $C_b^k(G;\R)$ by the space
$C^k(G;\R)$ of all functions $G^k\to \R$ and if the coboundary
operator is defined by the same formula, then we obtain the
cohomology $H^*(G;\R)$ of $G$. The inclusion
$C_b^k(G;\R)\hookrightarrow C^k(G;\R)$ induces a natural map
$H_{b}^k(G ; \R)\to H^k(G ; \R)$. When $k=2$, following
Grigorchuk~\cite{g}, let us denote the kernel of this map by
$H_{b,2}^2(G ; \R)$.

For a group $G$, let $PX(G)$ denote the space of pseudo characters
(pseudo homomorphisms) on $G$. That is, $PX(G)$ is the space of
all functions $f:G\to \R$ satisfying $|f(x)+f(y)-f(xy)|\leq C$ and
$f(x^n)=nf(x)$ for all $x,y\in G$ and for some $C$ depending on
$f$. Let $X(G)$ denote the space of all homomorphism $G\to \R$.
Grigorchuk proved that $H_{b,2}^2(G;\R)$ is isomorphic to
$PX(G)/X(G)$ as vector spaces.

The next lemma was proved in the proof of  Proposition $4.7$
in~\cite{g}.

 \begin{lem}
Let $G$ be a finitely generated group and let $H$ be a subgroup of
finite index in $G$. The map $\tau: PX(G)\to PX(H)$ induced by the
restriction is injective and the quotient space $PX(H)/\tau
(PX(G))$ is finite dimensional.
 \label{lem:subgroup}
 \end{lem}

 \begin{thm} [\cite{bo}]
Let $G$ and $F$ be two groups and let $\sigma :G\to F$ be an
epimorphism. Then $\sigma$ induces an injective linear map
$H_{b}^2(F ; \R)\to H_{b}^2(G ; \R)$.
 \label{th:quotientgp}
 \end{thm}

 \begin{thm} [\cite{mm}]
Suppose that $n\geq 2$. If $F_n$  is a free group of rank $n$,
then the dimension of the space $H_{b}^2(F_n;\R)$ is equal to the
cardinal of the continuum.
 \label{th:freegp}
 \end{thm}

 \begin{thm}
Let $G$ be a finitely presented group and let $H$ be a subgroup of
finite index in $G$. Suppose that there is a homomorphism from $H$
onto a free group $F_n$ of rank $n\geq 2$. Then the dimension of
the space $H_{b}^2(G;\R)$ is equal to the cardinal of the
continuum.
 \label{th:H_b^2}
 \end{thm}
 \begin{proof}
If $K$ is a finitely generated group, then it is countable. It
follows that the dimension of $C_b^k(K;\R )$, and hence that of
$H_b^k(K;\R)$, is at most the cardinal of the continuum for any
positive integer $k$.

Since $F_n$ is a quotient of $H$, it follows from
Theorems~\ref{th:quotientgp} and~\ref{th:freegp} that the
dimension of $H_{b}^2(H;\R)$ is the cardinal of the continuum.
Since $H$ is finitely presented, $H^2(H;\R)$ and $X(H)$ are finite
dimensional. It follows that the dimensions of $H_{b,2}^2(H;\R)$ and
$PX(H)$ are the cardinal of
the continuum. We conclude from Lemma~\ref{lem:subgroup} that
the dimension of $H_{b,2}^2(G;\R)$ is also the cardinal of the
continuum. Since $H_{b,2}^2(G;\R)$ is a subspace of
$H_{b}^2(G;\R)$, the theorem is follows.
 \end{proof}

 \begin{thm}
 Suppose that $g\leq 2$. If $g=0$, suppose, in addition, that
 $p+q\geq 4$. Then the dimension of $H_b^2(\mcg_{g,p}^q;\R)$ is
 equal to the cardinal of continuum.
 \label{th:g=0,1}
 \end{thm}

 \begin{proof}
The proof follows from Theorems~\ref{th:main} and~\ref{th:H_b^2}
and the fact that $\mcg_{g,p}^q$ is finitely presented.
 \end{proof}

 \begin{thm}
The dimension of the second bounded cohomology group
$H_b^2(\Delta_g;\R)$ of the hyperelliptic mapping
class group $\Delta_g$ is equal to the cardinal of
continuum.
 \label{th:hyp}
 \end{thm}

 \begin{proof}
The proof follows from Corollary~\ref{cor:main} and Theorem~\ref{th:H_b^2}
and the fact that $\Delta_g$ is finitely presented.
 \end{proof}

 \begin{rem}
 If $p+q\leq 3$ then the group $\mcg_{0,p}^q$ is
 either a finite group or a free abelian group. All these
 groups are amenable and amenable groups have trivial bounded
 cohomology.
 \end{rem}

\section{Automorphisms of cofinite subgroups of mapping class groups}

Let $S$ be a surface of genus $g$ with $p$ punctures.
It was shown in \cite{i1} and \cite{ko} that 
any isomorphism between two subgroups of finite index in 
the extended mapping class group of $S$
is the restriction of an automorphism of the extended mapping class group
provided that $p\geq 5$ if $g=0$, $p\geq 3$ if $g=1$ or $g\geq 2$. 
The extended mapping class group of $S$ is defined as the group of isotopy 
classes of all diffeomorphisms $S\to S$ including orientation reversing ones.
Since the mapping class group $\mcg_{g,p}$ is characteristic 
in the extended mapping class group, it follows that any 
isomorphism between two subgroups of finite index in
$\mcg_{g,p}$ is the restriction of an automorphism of $\mcg_{g,p}$
under the restrictions on $g$ and $p$ above. 
If $g=0$ and $p\leq 1$, then the mapping class group
is trivial. If $(g,p)=(0,2)$, then the mapping class group
is a cyclic group of order $2$. If $(g,p)=(0,3)$ then the mapping class group
is the symmetric group on three letters. Obviously, in
these cases any isomorphism between two subgroups of $\mcg_{g,p}$ is the
restriction of an automorphism of $\mcg_{g,p}$. We prove in this section that
this result does not hold if $(g,p)$ is equal to $(0,4)$, $(1,0)$ or $(1,1)$, 
leaving the case $(g,p)=(1,2)$ open.

 \begin{lem}
Let $F_n$ be a nonabelian free group of rank $n$. If $H$ is a proper subgroup of finite index
in $F_n$, then there exists an automorphism $\varphi : H\to H$
such that $\varphi$ is not the restriction of any endomorphism of $F_n$. 
 \label{lem:free}
 \end{lem}
\begin{proof}
Suppose that $F_n$ is generated by $\{ y, x_1,x_2,\ldots, x_{n-1}\}$.
Assume that the index of $H$ is $k\geq 2$, so that $H$ is a free group of rank $k(n-1)+1$.
Since any two such groups are isomorphic, we can assume without loss of generality
that $H$ is the subgroup generated by $\{ y^j x_i y^{-j}, y^k \ |\ 1\leq i\leq n-1,\  0\leq j\leq k-1\}$.

Define an automorphism $\sigma :H\to H$ by $\sigma (yx_1y^{-1})=yx_1y^{-1}x_1$ and the identity
on all other generators of $H$. The automorphism $\sigma$ does not extend to 
any endomorphism $F_n\to F_n$. Because if there is such an extension, then we 
conclude from $\tilde{\sigma}(y^k)=y^k$ that
$\tilde{\sigma}(y)=y$. Since $\tilde{\sigma}$ also fix all generators $x_i$ of $F_n$, 
it must be the identity. But $\sigma$ is not the identity.  
\end{proof}

 \begin{thm}
If $(g,p)$ is equal to $(0,4),(1,0)$ or $(1,1)$, then there exists a subgroup $\Gamma$ of
finite index mapping class group $\mcg_{g,p}$ and an automorphism 
$\varphi : \Gamma\to\Gamma$ such that $\varphi$ is not the restriction of any endomorphism
of $\mcg_{g,p}$.
 \label{thm:auto}
 \end{thm}
\begin{proof}
Suppose first that $g=1$. Note that in this case $\mcg_{1,0}$ and $\mcg_{1,1}$ 
are isomorphic to $SL(2,\Z)$. It is well known that the commutator subgroup of 
$SL(2,\Z)$, denoted by $F_2$, is a free group of rank $2$ and its index 
in $SL(2,\Z)$ is $12$. Let $\Gamma$ be any proper subgroup of finite index in $F_2$.
By Lemma~\ref{lem:free}, there exists an automorphism $\varphi : \Gamma\to\Gamma$ 
which is not the restriction of any endomorphism of $F_2$. 
Since any endomorphism $SL(2,\Z)$ induces an endomorphism
$F_2$, we are in this case.

Suppose now that $g=0$ and $p=4$. 
Let $S$ be a sphere with four punctures, say $P_1,P_2,P_3,P_4$.  For $i=1,2,3$ let $\alpha_i$ be
three disjoint  embedded arcs from $P_i$ to $P_{i+1}$. Let $a$ and $b$
denote the boundary component of a regular neighborhood of $\alpha_1$ and
$\alpha_2$, respectively. The pure mapping class group $\pmcg_{0,4}$ of
$S$ is a free group of rank two freely generated by the Dehn twists $t_a$ and $t_b$. 
Let $w_i$ denote the half twist about $\alpha_i$, so that $w_i$ interchanges $P_i$ and $P_{i+1}$,
$(w_1)^2=t_a$,  $(w_2)^2=t_b$, and $(w_3)^2$ is the right Dehn twist about the
boundary of a regular neighborhood of $\alpha_3$. 
Theorem~$4.5$ in \cite{birman} gives a presentation of $\mcg_{0,n}$ for all $n$.
It follows from this, in particular, that $w_1,w_2$ and $w_3$ generate
$\mcg_{0,4}$ and $H_1(\mcg_{0,4})$ is a cyclic group of order $6$
generated by the class of any $w_i$. Thus,
the classes of $t_a$ and  $t_b$ in $H_1(\mcg_{0,4})$ both have orders $3$.

We now define an automorphism $\varphi :\pmcg_{0,4}\to \pmcg_{0,4} $ by 
$\varphi (t_a)=t_a$ and $\varphi (t_b)=t_at_b$.
Suppose that there is an endomorphism $\tilde{\varphi}$ of 
$\mcg_{0,4}$ extending $\varphi$. Since $w_1$ and $w_2$ are conjugate,
so are $t_a=\tilde{\varphi} (w_1)^2$ and $t_at_b=\tilde{\varphi} (w_2)^2$. This
implies that the classes of $t_a$ and $t_at_b$ in $H_1(\mcg_{0,4})$
are equal. Therefore, $t_b$ represents $0$ in $H_1(\mcg_{0,4})$. 
By this contradiction, $\varphi$ cannot be extended
to an automorphism of $\mcg_{0,4}$.
\end{proof}

\end{document}